\documentclass[a4paper,11pt]{amsart}
\usepackage[latin1]{inputenc}
\usepackage[francais]{babel}
\usepackage{amsmath}
\usepackage{amsfonts}
\usepackage{amssymb}
\usepackage{hyperref}
\usepackage{multirow}
\usepackage{tabularx}

\newtheorem{theorem}{Théorème}[section]
\newtheorem{remark}{Remarque}[section]

\catcode`\@=11

\catcode`\;=\active
\def;{\relax\ifhmode\ifdim\lastskip>\z@\unskip\fi\kern.2em\fi\string;\ \ignorespaces}

\catcode`\:=\active
\def:{\relax\ifhmode\ifdim\lastskip>\z@\unskip\fi\penalty\@M\ \fi\string:\ \ignorespaces}

\catcode`\!=\active
\def!{\relax\ifhmode\ifdim\lastskip>\z@\unskip\fi\kern.2em\fi\string!\ \ignorespaces}

\catcode`\?=\active
\def?{\relax\ifhmode\ifdim\lastskip>\z@\unskip\fi\kern.2em\fi\string?\ \ignorespaces}

\newif\ifguill  \guilltrue
\catcode`\"=\active
\def"{\ifguill``\guillfalse\else''\guilltrue\fi}

\catcode`\@=12
\catcode`\à=\active \defà{\`a}%
\catcode`\À=\active \defÀ{\`A}%
\catcode`\â=\active \defâ{\^a}%
\catcode`\é=\active \defé{\'e}%
\catcode`\è=\active \defè{\`e}%
\catcode`\ê=\active \defê{\^e}%
\catcode`\ë=\active \defë{\"e}%
\catcode`\ï=\active \defï{\"\i}%
\catcode`\î=\active \defî{\^\i}%
\catcode`\ô=\active \defô{\^o}%
\catcode`\ö=\active \defö{\"o}%
\catcode`\ù=\active \defù{\`u}%
\catcode`\û=\active \defû{\^u}%
\catcode`\ü=\active \defü{\"u}%
\catcode`\ç=\active \defç{\c c}%

\newcommand{\N}{\mathbb{N}}
\newcommand{\Z}{\mathbb{Z}}
\newcommand{\Q}{\mathbb{Q}}
\newcommand{\R}{\mathbb{R}}
\newcommand{\C}{\mathbb{C}}
\newcommand{\s}{\mathbb{S}}
\newcommand{\be}{\begin{enumerate}}
\newcommand{\ee}{\end{enumerate}}
\newcommand{\soe}{\geq} 
\newcommand{\ioe}{\leq} 
\newcommand{\tq}{\text{ tel que }}
\newcommand{\vers}{\longrightarrow}
\newcommand{\app}[5]{\begin{array}[t]{lrll}
					#1&#2&\longrightarrow & #3 \\
					&#4&\longmapsto & #5
			\end{array}} 
\newcommand{\qlq}{\forall\;}
\newcommand{\impliq}{\Longrightarrow}
\newcommand{\ssit}{\textit{ si et seulement si }}
\newcommand{\ssi}{\Longleftrightarrow}
\newcommand{\ds}{\displaystyle}
\begin{document}

\title[La Conjecture H]{La conjecture H: Une minoration de la dimension cohomologique pour un espace elliptique
\\ \today
}
\vskip 10mm
\author{M.R Hilali}
\address{D\'epartement de Math\'ematiques\\
         Faculté des sciences Ain Chok\\
         Universit\'e Hassan II,  Route d'El Jadida\\
         Casablanca\\
         Maroc}
\email{rhilali@hotmail.com}

\author{M.I Mamouni}
\address{Classes préparatoires aux grandes écoles d'ingénieurs\\
         Lycée Med V\\
         Avenue 2 Mars\\
		Casablanca\\
        Maroc}
\email{myismail1@menara.ma}

\vskip 20mm
\begin{abstract}
Le but de cet article est d'améliorer les conditions suffisantes, déjà établi par le premier auteur pour que la somme des nombres de Betti, d'un CW-complexe 1-connexe fini et rationnel, soit supérieur à la dimension de son $\Q$-espace vectoriel d'homotopie, qu'on présentera dans deux aspects, celui algébrique et un autre géométrique.
\end{abstract}



\maketitle

\tableofcontents

\newpage
\section{Introduction}
\subsection{Les CW-complexes finis  1-connexes} Les CW complexe finis  \\1-connexes, notés ci-après $X$,  se divisent en deux classes:

\begin{itemize}
\item Les CW complexes  elliptiques  tels que $$\dim \pi_{\ast}(X) \otimes \mathbb Q < \infty\,.$$

\item Les CW complexes hyperboliques   tels que $$ \dim \pi_{\ast}(X) \otimes \mathbb Q = \infty \,.$$
\end{itemize}
\subsection{Les elliptiques.} La cohomologie rationnelle $H^\ast (X;\mathbb Q)$ d'un espace elliptique $X$  est à dualité de Poincaré et sa caractéristique d'Euler, $\chi $ vérifie  $\chi \geq 0$.
\vskip 3mm
Bien que les espaces elliptiques jouissent de propriétés très spécifiques  (nous en verrons d'autres au cours de cette exposé)  ce sont ceux que l'on rencontre le plus couramment en géométrie différentielle (Groupes de Lie, espaces homogènes,...).
\subsection{Quelques conjectures dans le cas elliptique.} J. Moore  a conjecturé : 

\textbf{Conjecture M}. \textit{Le groupe abélien gradué $\pi_\ast(X)$ possède un $p$-exposant ssi $X$ est  CW complexe  elliptique.}

\vskip 3mm
Rappelons qu'un groupe gradué vérifie une propriété si chacune de ses composantes vérifie cette propriété et qu'un groupe abélien admet un $p$-exposant ($p$ nombre premier) s'il existe un entier $n$ tel que $p^n$ annule la $p$-torsion du groupe. Cette conjecture a été établie dans de nombreux cas particuliers par P. Sélick, Thériault, Stelzer,... .
\vskip 3mm
R. Bott a conjecturé que

\textbf{Conjecture B}. \textit{Toute variété riemannienne compacte 1-connexe sans bord dont la courbure sectionnelle est toujours $\geq 0$ est un CW complexe elliptique.}

Cette conjecture a été partiellement démontrée par G. Paternain en 1992.

\vspace{3mm}
Dans cet article nous sommes intéressés par la conjecture H, émise par le premier auteur en (1990) :

\vspace{3mm}
\textbf{Conjecture H}. \textit{Pour tout CW complexe elliptique, 1-connexe $X$, on a:}
$$\dim \pi_\ast(X) \otimes \mathbb Q \leq \dim H^\ast(X; \mathbb Q)\,.
$$
Cette conjecture été démontrée par le premier auteur \cite{Hi} pour les CW complexes elliptiques tels que $\chi >0$,  (cas pur).
\subsection{Contexte.}
La conjecture H, entre dans le cadre d'autres conjectures qui proposent une minoration de la dimension cohomologique, on citera en particulier celle du rang torique due à S. Halperin en 1986.
\vskip 3mm

\textbf{Conjecture du rang torique (CRT).} \textit{Si $ X$ est un espace 1-connexe et raisonnable, alors}
$$\dim H^*(X,\Q)\soe 2^{rk_0(X)}$$

On rappelle qu'un espace $X$ est dit raisonnable s'il vérifie les propriétés suivantes:
\begin{itemize}
\item $X$ est connexe, de Hausdorff, compact ou paracompact.
\item $\dim H^*(X,\Q)<\infty$.
\item Pour tout $x\in X$, on a $\ds\lim_{\vers} H^*(U,\Q)=0$ où la limite est prise sur les voisinages ouverts $U$ de $x$.
\end{itemize}

Cette conjecture a été résolue dans des cas particuliers par Halperin, Allday, Puppe, Hilali,...

\subsection{La conjecture H, dans le cas non elliptique.}
Toutes les situations sont possibles, on les illustrera ci-dessous avec des exemples:
\begin{itemize}
\item $\dim H^*(X,\Q)=+\infty,\dim \pi_\ast(X) \otimes \mathbb Q =+\infty$.

Prendre $X=\s^2\times\s^3\times\cdots$
\item $\dim H^*(X,\Q)=+\infty,\dim \pi_\ast(X) \otimes \mathbb Q <+\infty$.

Prendre $X=\C {\rm P}^\infty$, dans ce cas

$H^*(X,\Q)=\Q[a]\tq |a| \text{ pair}$ et $\pi_*(X) \otimes \mathbb Q=K(\Z,2)$

\item $\dim H^*(X,\Q)<+\infty,\dim \pi_\ast(X) \otimes \mathbb Q =+\infty$.

Prendre $X=\s^3\vee \s^3$, dans ce cas 

$H^0(X,\Q)\equiv\Q,\; H^3(X,\Q)\equiv \Q\oplus\Q$.

$\pi_{*+1}(\s^3\vee \s^3)=\pi_{*}(\Omega(\s^3\vee \s^3))={\mathbb{L}}(a,b)$ avec $ |a|=|b|=2$.
\end{itemize}
\subsection{Notre contribution.}
Nous démontrerons la conjecture H  dans quelques cas où  $\chi =0$.

\vskip 3mm
Rappelons que d'après \cite{FHT}-[Proposition 32.16], 
 $$
 \left\{
 \begin{array}{ll}
 \chi >0&\Longleftrightarrow  \dim \pi_{impair} \otimes \mathbb Q = \dim  \pi_{pair} \otimes \mathbb Q\\
 \chi =0&\Longleftrightarrow  \dim \pi_{impair} \otimes \mathbb Q > \dim  \pi_{pair} \otimes \mathbb Q\\
 \end{array}
 \right.
 $$ 
et donc que le cas restant, $\chi =0$,  est le plus vaste.
\vskip 3mm

Plus précisément, nous suivrons le plan suivant:
\begin{itemize}
\item Quelques points de la théorie des modèles minimaux de Sullivan.
\item Enoncé de la conjecture en termes de modèles minimaux.
\item Le cas hyper-elliptique, sous condition. (Enoncé et preuve du théorème).
\item Le cas elliptique, sous condition. (Enoncé et preuve du théorème).
\item Le rang torique. (Enoncé et preuve du théorème).
\item Le cas  symplectique, sous condition. (Enoncé et preuve du théorème).
\item Le cas co-symplectique, sans condition. (Enoncé et preuve du théorème).
\end{itemize}
\section{Résolution du problème.}
\subsection{Motivation}
La conjecture (H) est vérifiée pour les H-espaces, car l'image du morphisme d'Hurewicz, $hur_X:\pi_*(X)\otimes \Q\vers H_*(X,\Q)$ engendre l'algèbre de Pontrayagin, $H_*(X,\Q)$ et car $H_*(X,\Q)\cong H^*(X,\Q)$ par dualité.
\subsection{Algèbrisation de la conjecture.} D.Sullivan \cite{Su} a construit pour tout espace topologique 1-connexe, $X$, vérifiant  $\dim H^i (X;\mathbb Q)<\infty $ pour tout $i$, un modèle minimal $(A_{PL}(X),d)=(\bigwedge V,d)$, unique à isomorphisme prés quand $H^0(X,\Q)=\Q$, tel que \cite{Ha2}:
$$\begin{array}{rcl}
V^n&\cong&\pi_n(X)\otimes\Q\\
H^n(\bigwedge V,d)&\cong&H^n(X,\Q)\\
\end{array}$$
donc 
$$\begin{array}{rcl}
V&\cong&\pi_*(X)\otimes\Q\\
H^*(\bigwedge V,d)&\cong&H^*(X,\Q)\\
\end{array}$$
Ce qui nous permet d'énoncer  la conjecture H de la façon suivante:

\textbf{Version algèbrique.} Si $(\bigwedge V,d)$ désigne un modèle minimal tel que:\\
$\begin{array}{ll}
i)&\dim(V)<\infty\\
ii)&\dim H^*\left(\bigwedge V,d\right)<\infty
\end{array}$
\\Alors :
$$
 \dim(V)\ioe \dim H^*\left(\bigwedge V,d\right)
$$
\subsection{Les théorèmes.}
La conjecture (H) est vérifiée dans les cas suivants:

\vskip 3mm
\textbf{Théorème A}. $X$ \textit{est hyper-elliptique vérifiant:}
$$\dim\left(\pi_{\text{pair}}\otimes \Q\right)\soe \frac 12\left(1+\sqrt{-12\chi_\pi(X)-15}\right)$$

\vskip 3mm
\textbf{Corollaire A.1}. $X$ \textit{est hyperelliptique vérifiant:}
\begin{equation}\label{CS2_ConjH}
\chi_\pi(X)\in\{0,-1,-2\}
\end{equation}

\vskip 3mm
\textbf{Théorème B.} $X$ \textit{est elliptique 1-connexe vérifiant:}

\begin{equation}]\label{CS3_ConjH}
fd(X)\ioe 10
\end{equation}

Où $fd(X)$ appelée dimension formelle de $X$ est définie par la relation suivante:
$$fd(X)=\max\{n\in\N\tq H^k(X,\Q)\neq 0\}$$

\vskip 3mm
\textbf{Théorème C.}  $X$ \textit{est un espace rationnel 1-connexe, elliptique tel que: 
$\pi_{pair}(X)\otimes \Q=0$, de modèle $\bigwedge V=\bigwedge\{y_1,\ldots ,y_n\}$ avec $|y_i|$ impair vérifiant:}
$$\qlq i\in\{1,\ldots n\}\quad \dim (ker \delta_i)>\dim (Im \delta_i)$$

Où $\app{\delta_i:}{A_{i-1}}{A_{i-1}}\beta{\beta\alpha_i}$ avec $A_i=H^*(\wedge\{y_1,\ldots ,y_i\},d)$ et $\alpha_i=[dy_i]$.

\vskip 3mm
\textbf{Corollaire C.1.}  \textit{$X$
elliptique, 1-connexe de mod\`{e}le minimal $(\bigwedge V,d)=(\bigwedge
\{y_{1},y_{2},...,y_{n}\},d)$  v\'{e}rifiant: }

\begin{itemize}
\item $\pi _{2k}(X)\otimes \mathbb{Q}=0$, $\forall k\in \N$ 

\item $\alpha_i^2 = [dy_i]^2 = 0$, $\forall i\in \{3,...,n\}$ 

\item $\forall i\in \{3,...,n\}$,\textit{ il existe} $%
\gamma_{1i},\gamma_{2i}\in A_i^+$, \textit{tels que}: $\alpha_i =
\gamma_{1i}\gamma_{2i}$, \textit{et} $\gamma_{1i}^2=0$. 
\end{itemize}

\vskip 3mm
\textbf{Théorème D.}  $X$ \textit{est hyperelliptique  vérifiant:}
$$rk_0(X)=-\chi_\pi(X)-i\tq i\in\{0,1,2\}$$

Où $rk_0(X)$ appelé rang torique de $X$ est défini par la relation suivante:
$$rk_0(X)=\max\{n\in\N\tq {\mathbb T}^n=\left( {\mathbb S}^1\right)^n\text{ agit presque librement sur }X\}$$

c'est à dire: $\qlq x\in X, G_x=\{g\in {\mathbb T}^n\tq g.x=x\}$ est fini.

\vskip 3mm
\textbf{Théorème E.} $X$ \textit{est elliptique 1-connexe vérifiant:}
$$fd(X)-rk_0(X)\ioe 6$$

$codim(X)=fd(X)-rk_0(X)$  s'appelle codimension de $X$.

\vskip 3mm
\textbf{Théorème F.} $X$ \textit{ est 1-connexe, symplectique.}

\vskip 3mm
\textbf{Théorème G.} $X$\textit{ est  co-symplectique}.

\subsection{Les outils.}
\subsubsection{Définitions, notations et vocabulaire}.

\vskip 3mm
\textbf{Modèles de Sullivan.}
\begin{itemize}
\item $(\bigwedge V,d)$ est dit modèle minimal elliptique, \ssit :
\be
\item $\dim V<\infty$.
\item $\dim H^*\left(\bigwedge V,d)\right)<\infty$.
\item Il existe $\{|x_1|\ioe \cdots\ioe |x_n|\}$ base de $V$ telle que
$$
\begin{array}[t]{ll}
i)&dx_1=0\\
ii)&dx_i\in\bigwedge^{\soe 2}\{x_1,\ldots x_{i-1}\},\quad\qlq 1\ioe i\ioe n
\end{array}
$$
\ee
\item Le modèle est dit pur si de plus: $$d(V^{\text{pair}})=0\text{ et }d(V^{\text{impair}})\subset \bigwedge V^{\text{pair}}$$
\item Le modèle est dit hyperelliptique si de plus 
\begin{equation}
\label{def_hyp}
\begin{array}[t]{l}
d(V^{\text{pair}})=0\\
d(V^{\text{impair}})\subset \bigwedge^+V^{\text{pair}}\otimes\bigwedge V^{\text{impair}}
\end{array}
\end{equation}

\end{itemize}

\vskip 3mm
\textbf{Invariants d'Euler-Poincaré.}
Dans le cas d'un espace elliptique, on définit les deux invariants suivants:
\begin{itemize}
\item \textit{Invariant cohomologique:}
$$\begin{array}{ll}
\chi_c(X)&=\ds\sum_{k\soe 0}(-1)^k\dim H^k(X,\Q)\\
&=\dim \left(H^{\text{pair}}(X,\Q)\right)-\dim\left(H^{\text{impair}}(X,\Q)\right)
\end{array}$$

\item \textit{Invariant homotopique:}
$$
\begin{array}{ll}
\chi_\pi(X)&=\ds\sum_{k\soe 0}(-1)^k\dim\pi_k(X)\otimes\Q\\
&=\dim(\pi_{\text{pair}}(X)\otimes\Q)-\dim(\pi_{\text{impair}}(X)\otimes\Q)\\
&=\dim(V^{\text{pair}})-\dim(V^{\text{impair}})
\end{array}$$
\end{itemize}

\vskip 3mm
\textbf{H-espaces} \cite{Th}, \cite{Hat-02} Un $H$-espace  est un espace topologique, $X$, muni d'une application continue $\mu:X\times X\vers X\tq \mu\circ i_\varepsilon =id_X$  où $i_\varepsilon: X\vers X\times X, \varepsilon=1,2$ désigne l'une des inclusions naturelles.
\vskip 3mm

Les H-espaces abondent en géométrie et en topologie, comme exemple on peut citer:

\begin{itemize}
\item Les groupes topologiques, en particulier les groupes de Lie.

\item Les sphères ${\mathbb S}^0$, ${\mathbb S}^1$, (des complexes), ${\mathbb S}^3$, (des quaternions), ${\mathbb S}^7$, (des octanions). 
\\Adams \cite{Ad} a démontré que ce sont les seules sphères, $H$-espaces.

\item ${\mathbb R}{\rm P}^1={\mathbb S}^1\diagup \pm 1, {\mathbb R}{\rm P}^3={\mathbb S}^3\diagup \pm 1,{\mathbb R}{\rm P}^7={\mathbb S}^7\diagup \pm 1$.
\\Dans le cas général ${\mathbb R}{\rm P}^n$ est un $H$-espace \ssit $n+1$ est une puissance de 2.

\item L'espace de lacets $\Omega X$, d'un  espace pointé $X$.

\item Les espaces d'Eilenberg-MacLane, $K(G,n)$, tel que $n\soe  1$ et $G$ abélien, puisque $K(G,n)=\Omega K(G,n+1)$.

\item $\C {\rm P}^\infty$.

\item $J(X)$: The "James reduced product" associé à un espace topologique pointé, $(X,*)$, défini par la relation: $$J(X)=\ds\coprod_{k\soe 1}X^k\diagup (x_1,\ldots,x_i,\ldots,x_k)\sim (x_1,\ldots,\widehat{x_i},\ldots,x_k),\text{ si }x_i=*$$

\item $SP(X)$: The "infinite symmetric product" associé à un espace topologique, $X$, défini par la relation: $$SP(X)=\ds\coprod_{k\soe 1}X^k\diagup (x_1,\ldots,x_k)\sim (x_{\sigma(1)},\ldots,x_{\sigma(k)})$$
\end{itemize}

\vskip 3mm
\textbf{Espaces symplectiques} \cite{La} Une variété symplectique  est une variété différentielle $M$, munie d'une forme différentielle de degré deux $\omega$ fermée et non dégénérée, appelée forme symplectique.

L'étude des variétés symplectiques relève de la topologie symplectique. Les variété symplectiques apparaissent lors de l'étude des formulations abstraites de la mécanique classique (Hamiltonienne) et analytique (quantique), liées au fibrés cotangents des variétés, notamment dans la description hamiltonienne de la mécanique, où les configurations d'un système forment une variété dont le fibré cotangent décrit l'espace des phases du système.

Les espaces de Kahler sont des variétés symplectiques, qui  sont des variétés de Poisson, qui à leur tour sont des variétés de Jacobi.

\vskip 3mm
\textbf{Espaces cosymplectiques}. Une structure cosymplectique sur une variété $M$ de dimension $2k+1$ est la donnée d'une 1-forme fermée $\theta$ et une 2-forme fermée $\omega$ telles que 
$\theta\wedge\omega^k$ est une forme volume sur $M$, où $\omega^k$ désigne le produit de $k$ copies de $\omega$.

\subsubsection{Principaux résultats utilisés.}
\begin{theorem}\label{thm_AH}
C. Allday \& S. Halperin, \cite{AH} 

Si $X$ est elliptique 1-connexe et raisonnable, alors:
\begin{itemize}
\item $\chi_\pi(X)\ioe -rk_0(X).$
\item $\chi_\pi(X)=-rk_0(X)\impliq X$ pur.
\end{itemize}
\end{theorem}
\begin{theorem} \label{thm_FH} J. Friedlander \& S. Halperin, \cite{FH} 

Si $X$ espace elliptique, 1-connexe de modèle minimal de Sullivan, $(\bigwedge V,d)$, alors 
$$\dim V\ioe fd(X)$$

Et on peut trouver  $\{x_1,\ldots,x_n\}$ base de $V^{\text{pair}}$ et  $\{y_1,\ldots,y_{n+p}\}$ base de $V^{\text{impair}}$ telles que:
$$
\begin{array}{l}
|x_1|\ioe\cdots\ioe |x_n|,\;\qlq 1\ioe i\ioe n\\
\\
|y_i|\soe 2|x_i|-1,\;\qlq 1\ioe i\ioe n\\
\ds\sum_{i=1}^n|x_i|\ioe fd(X)\\
\ds\sum_{i=1}^{n+p}|y_i|\ioe 2fd(X)-1\\
\ds\sum_{i=1}^{n+p}|y_i|-\ds\sum_{i=1}^n(|x_i|-1)=fd(X)
\end{array}
$$
\end{theorem}
\begin{theorem} \label{thm_Ha}S. Halperin.\cite{Ha2}

Si $X$ espace elliptique, 1-connexe alors:
\begin{itemize}
\item $\chi_\pi(X)\ioe 0,\;\chi_c(X)\soe 0$
\item $
\begin{array}[t]{lll}
\chi_\pi(X)< 0&\ssi &\chi_c(X)= 0
\end{array}$
\end{itemize}
\end{theorem}
On prendra dans la suite $\{x_1,\ldots,x_n\}$ comme base de $V^{\text{pair}}$ et 
$\{y_1,\ldots,y_{n+p}\}$ comme base de $V^{\text{impair}}$, avec $p\soe 0$ et $\chi_\pi=-p$.

En particulier on peut en conclure  que:
\begin{itemize}
\item $\chi_c(X)=0\impliq\dim H^{*}(X,\Q)=2\dim H^{\text{pair}}(X,\Q)$.
\item $\chi_c(X)\neq 0\impliq\dim H^{*}(X,\Q)=\dim H^{\text{pair}}(X,\Q)$.
\end{itemize}
\begin{theorem} \label{thm_Hi}M.R. Hilali, \cite{Hi}

Soit $X$ un CW-complexe 1-connexe,
v\'erifiant l'une des deux conditions suivantes: 

\begin{itemize}
\item $X$ est un espace hyperelliptique.

\item $X$ est de dualit\'e de Poincar\'e,  tel que ${\rm codim}(X) \leq 6$. 
\end{itemize}
Alors: 
$$\dim H^{\star }(X,\mathbb{Q})\geq 2^{rk_{0}(X)}.$$

\end{theorem}
\begin{theorem} S. Halperin \& G. Levin, \cite{HL}

Si $X$ est elliptique 1-connexe et pur, alors: $$\chi_\pi(X)=-rk_0(X).$$
\end{theorem}

\section{Les démonstrations.}
\subsection{Cas hyper-elliptique} (sous condition)

\vskip 3mm

\textbf{Preuve du théorème A.} D'après (\ref{def_hyp}) on a

$$
\begin{array}{ll}
dx_i=0 &\qlq i\in\{1,\ldots,n\}\\
dy_j=P_j+w_j&\qlq j\in\{1,\ldots,n+p\}\tq P_j\in\bigwedge^+V^{\text{pair}},w_j\in\bigwedge^+V^{\text{pair}}\otimes \bigwedge^+ V^{\text{impair}}
\end{array}
$$

Notons $W_1, W_2$ les sous-espaces vectoriels de $H^{\text{pair}}(\bigwedge V,d)$ engendrés respectivement par $\left([x_i]\right)_{1\ioe i\ioe n}$ et $\left([x_ix_j]\right)_{1\ioe i\ioe j\ioe n}$. Grâce à la minimalité du modèle on a  $H^0(\wedge V,d)\oplus W_1\oplus W_2$ est une somme directe dans $H^{\text{pair}}(\wedge V,d)$, et  $\left([x_i]\right)_{1\ioe i\ioe n}$ libre, donc
$$\begin{array}{l}
H^0(\wedge V,d)\oplus W_1\oplus W_2\subset H^{\text{pair}}(\wedge V,d)\\
\dim H^0(\wedge V,d)=1,\dim W_1=\dfrac{n(n+1)}2
\end{array}$$

D'autre part 
$$\begin{array}{l}
W_2\oplus(\bigwedge^2V^{\text{pair}}\cap dV^{\text{impair}})=\bigwedge^2V^{\text{pair}}\\
\dim \bigwedge^2V^{\text{pair}}=\dfrac{n(n+1)}2\\
\dim\bigwedge^2V^{\text{pair}}\cap dV^{\text{impair}}\ioe \dim dV^{\text{impair}}=n+p
\end{array}$$
car $\{dy_1,\ldots,dy_{n+p}\}$ engendre $dV^{\text{impair}}$

Supposons que $X$ n'est pas pur (le cas contraire a été déja résolu), donc $\exists j\in\{1,\ldots,n+p\}\tq w_j\neq 0$, d'où $dy_j\notin \wedge ^2V^{\text{pair}}$, donc $dy_j\notin \bigwedge^2V^{\text{pair}}\cap dV^{\text{impair}}$, et par suite  
\begin{equation}\label{ineq}
\begin{array}{l}
\dim V_1\ioe n+p-1
\\
\dim W_2\soe \dfrac{n(n+1)}2-n-p+1\\
\dim H^{\text{pair}}(\bigwedge V,d)\soe \dfrac{n(n+1)}2-p+2\\
\end{array}
\end{equation}

Comme $X$ n'est pas pur, alors $\chi_\pi(X)<-rk_0(X)\ioe 0$ (théorème \ref{thm_AH}), donc $\chi_c(X)=0$ (théorème \ref{thm_Ha}), donc $\dim H^*(\bigwedge V,d)=2\dim H^{\text{pair}}(\bigwedge V,d)$, or $\dim V=2n+p$, (\ref{ineq}) devient alors
\begin{equation}\label{ineq_trinom}
\dim H^*(\bigwedge V,d)-\dim V\soe n^2-n-3p+4
\end{equation}

L'étude du signe du trinôme $P(n,p)=n^2-n-3p+4$, permet de conclure que la conjecture (H) est vraie pour 
\begin{equation}\label{CS_ConjH}
\begin{array}{lll}
p\soe 2&\text{ et }&n\soe\dfrac{1+\sqrt{12p-15}}2\\
p=1&\text{ et }&\qlq n\in\N\\
\end{array}
\end{equation}

Signalons enfin que $p=0$ est le cas pur, et que 
\\$n=\dim V^{\text{pair}}=\dim(\pi_{\text{pair}}\otimes\Q), \;p=-\chi_\pi(X)$. Ce qui termine notre démonstration \hfill $\Box$.

\vskip 3mm

\textbf{Preuve du Corollaire A.1.} D'aprés (\ref{ineq_trinom}) il suffit d'étudier le signe de $P(n,p)=n^2-n-3p+4$, on distingue les cas suivants 

\textbf{Premiers cas: }  $\chi _{\pi }(X)=0$, cas pur (déjà résolu).

\textbf{Deuxième cas: } $\chi _{\pi }(X)=-p=-1$\\On a: $P(n,1)=n^2-n+1\soe 0,\;\qlq n\in\N$.

\textbf{Troisième cas: }$\chi _{\pi }(X)=-p=-2$\\On a: $P(n,2)=n^2-n-2\soe 0,\;\qlq n\soe 2$. Pour $n\in\{0,1\}$, on a $\dim H^*(\bigwedge V,d)\soe 2(1+\dim W_1)=2+2n=\dim V$

\hfill $\Box$.

\subsection{Cas elliptique} (sous condition)

\vskip 3mm
\textbf{Preuve du Théorème B.}
Soit $\{x_1,\cdots,x_n\}$ une base de $V^{\text{pair}}$ et $\{y_1,\cdots,y_{n+p}\}$ une base de $V^{\text{impair}}$, d'après théorème \ref{thm_FH}, on a les tableaux récapitulatifs suivants:

\begin{center}
\begin{tabular}{|l|l|l|}
\hline
$fd(X)$ & $(|x_1|,\cdots,|x_n|)$ & $(|y_1|,\cdots,|y_{n+p}|)$\\
\hline
2&(2)&(3)\\
\hline
3& 0&(3)\\
\hline
\multirow{3}{2cm}{4}
 & (2)   & (5) \\  \cline{2-3}
 & (2,2)   & (3,3) \\  \cline{2-3}
 & (4)     & (7) \\
\hline
\multirow{2}{2cm}{5}
 & 0     & (5) \\  \cline{2-3}
 & (2)       & (3,3) \\  \cline{2-3}
\hline
\multirow{5}{2cm}{6}
 & 0     & (3,3) \\  \cline{2-3}
 & (2)       & (7) \\  \cline{2-3}
& (2,2)     & (3,5) \\  \cline{2-3}
 & (2,4)       & (3,7) \\  \cline{2-3}
 & (2,2,2)       & (3,3,3) \\  \cline{2-3}
\hline
\multirow{4}{2cm}{7}
 & 0     & (7) \\  \cline{2-3}
 & (2)       & (3,5) \\  \cline{2-3}
& (4)     & (7.3) \\  \cline{2-3}
 & (2,2)       & (3,3,3) \\  \cline{2-3}
\hline
\multirow{7}{2cm}{8}
 & 0     & (3,5) \\  \cline{2-3}
 & (2)       & (3,3,3) \\  \cline{2-3}
& (2,2,2)     & (3,3,5) \\  \cline{2-3}
 & (2,4)       & (5,7) \\  \cline{2-3}
 & (2,2,2,2)       & (3,3,3,3) \\  \cline{2-3}
 & (4,4)       & (7,7) \\  \cline{2-3}
 & (6)       & (13) \\  \cline{2-3}
\hline
\multirow{8}{2cm}{9}
 & 0     & (9) \\  \cline{2-3}
 & (2)       & (3,7) \\  \cline{2-3}
& (2)     & (5,5) \\  \cline{2-3}
 & (2,2)       & (3,3,5) \\  \cline{2-3}
 & (2,2,2)       & (3,3,3,3) \\  \cline{2-3}
 & (2,4)       & (3,7,3) \\  \cline{2-3}
 & (4)       & (7,5) \\  \cline{2-3}
 & (6)       & (11,3) \\  \cline{2-3}
\hline
\end{tabular}
\end{center}
On vérifie  à la main, cas par cas que l'espace est pur, donc vérifie la conjecture (H).
\begin{center}
\begin{tabular}{|l|l|l|}
\hline
$fd(X)$ & $(|x_1|,\cdots,|x_n|)$ & $(|y_1|,\cdots,|y_{n+p}|)$\\
\hline
\multirow{13}{2cm}{10}
 & 0     & (5,5) \\  \cline{2-3}
 & 0       & (3,7) \\  \cline{2-3}
& (2)     & (3,3,5) \\  \cline{2-3}
 & (2,2)       & (3,3,3,3) \\  \cline{2-3}
 & (2,2,2)       & (3,5,5) \\  \cline{2-3}
 & (2,2,2)       & (3,3,7) \\  \cline{2-3}
 & (2,2,2,2)       & (3,3,3,5) \\  \cline{2-3}
 & (2,2,2,2,2)       & (3,3,3,3,3) \\  \cline{2-3}
 & (2,4)       & (7,7) \\  \cline{2-3}
 & (2,6)       & (5,11) \\  \cline{2-3}
 & (2,4,4)       & (3,7,7) \\  \cline{2-3}
 & (4)       & (7,3,3) \\  \cline{2-3}
 & (4,6)       & (7,11) \\  \cline{2-3}
\hline
\end{tabular}
\end{center}

On vérifie  à la main, cas par cas que l'espace est  hyperelliptique avec $\chi_\pi\in\{0,-1,-2\}$,  donc vérifie la conjecture (H) d'après \ref{CS2_ConjH} \hfill $\Box$.
\vskip 3mm
\textbf{Preuve du Théorème C}.
On va raisonner par récurrence, il est clair que c'est vrai pour $n=1$.\\
Supposons le résultat vrai pour $n-1$, donc $\dim H^*(\bigwedge W,d)\soe \dim W=n-1$ où $\bigwedge W=\bigwedge\{y_1,\ldots ,y_{n-1}\}$.\\
Considérons la suite courte $$0\vers (\bigwedge W,d)\vers (\bigwedge V=\bigwedge W\otimes \bigwedge y_n,d)\vers (\bigwedge W,d)\vers 0$$
Elle induit en cohomologie la suite exacte longue de connectant \cite{Su}
$$\app{\delta_n:}{A_{n-1}=(\bigwedge W,d)}{A_{n-1}=(\bigwedge W,d)}\beta{\beta\alpha_n}\qquad \alpha_n=[dy_n]$$
Qui nous induit à son tour la suite exacte courte
$$0\vers coker\;\delta_n\vers H^*(\bigwedge V,d)\vers ker\;\delta_n\vers 0$$ 

Ainsi $\begin{array}[t]{ll}
\dim H^*(\bigwedge V,d)&=\dim\;ker\;\delta_n+\dim\;coker\;\delta_n\\
&=\dim\;ker\;\delta_n+\dim H^*(\bigwedge V,d)-\dim\;Im\;\delta_n\\
&=2\dim\;ker\;\delta_n\\
&>\dim\;ker\;\delta_n+\dim\;Im\;\delta_n\\
&=\dim H^*(\bigwedge V,d)\\
&\soe \dim W=n-1
\end{array}$\\
Donc $\dim H^*(\bigwedge V,d)\soe n=\dim V$ \hfill $\Box$.

\vskip 3mm

\textbf{Preuve du Corollaire C.1}
On va montrer que $\forall i,1\geq i\geq n$, on a: $$\dim (ker\delta _{i})>\dim (Im\delta _{i})$$

En effet d'une part si $\alpha _{i}=0$, alors $\dim (Im\delta _{i})=0<\dim (ker\delta _{i})$, d'autre part si $\alpha _{i}\neq 0$ alors $\delta
_{i}\circ \delta _{i}=0$ et donc $Im\delta _{i}\subset ker\delta
_{i}$, montrons que cette inclusion est stricte. 

Supposons par l'absurde que $Im(\delta _{i})=ker(\delta _{i})$ alors $\forall
\beta \in Im(\delta _{i})\setminus\{0\},|\beta |\geq |\alpha _{i}|$. Or d'apr\`{e}s les hypoth\`{e}ses $\alpha _{i}=\gamma _{1i}\gamma
_{2i},\gamma _{1i}^{2}=0$ et donc $\alpha
_{i}\gamma _{1i}=\gamma _{2i}\gamma _{1i}^{2}=0$, ce qui implique
que $\gamma _{1i}\in ker(\delta _{i})$, ceci est impossible car $|\gamma _{1i}|<|\alpha _{i}|.$ \hfill $\Box$

\subsection{Cas du rang torique.} (sous condition)

\vskip 3mm

\textbf{Preuve du théorème D.}
Gardons les m\^{e}mes notations ci-dessus. Nous allons distinguer
trois cas:

\textbf{Premier cas: }$rk_{0}(X)=-\chi _{\pi }(X)=p.\; X$.
\\Cas pur d'après le théorème \ref{thm_AH}

\textbf{Deuxi\`{e}me cas:} $rk_{0}(X)=-\chi _{\pi }(X)-1=p-1$. On sait que:
$$\begin{array}{ll}
\mathit{\dim H}^{\star }\mathit{(}\bigwedge \mathit{V,d)\geq 2}^{p-1}&\text{d'après   théorème }\ref{thm_Hi}\\
\mathit{\dim H}^{\star }\mathit{(}\bigwedge \mathit{V,d)\geq
\dim V,\;\;\forall n\geq }\frac{1+\sqrt{12p-15}}{2}&\text{d'après }(\ref{CS_ConjH})
\end{array}$$

Par cons\'{e}quent 
\begin{equation}\label{CS1_ConjH}
\dim H^{\star }(\bigwedge V,d)\geq \dim V,\;\forall n\in \left[0,\frac{2^{p-1}-p}{2}\right]\cup \left[ \frac{1+\sqrt{12p-15}}{2},+\infty \right[
\end{equation}

Cherchons d'abord pour quelles valeurs de $p$ on a:
\begin{equation*}
\hspace{-3cm}\frac{2^{p-1}-p}{2}\geq \frac{1+\sqrt{12p-15}}{2}
\end{equation*}

et donc cherchons $p\in \N$ tels que:
$$A_{p}=2^{p}(2^{p-2}-p-1)+p^{2}-10p+16\geq 0$$

$\forall p\geq 8$ on a $A_{p}>0$, car  $p^{2}-10p+16\geq 0$ et $2^{p-2}-p-1>0$. On v\'{e}rifie facilement que $A_{p}>0,$ si $p\in \{5,6,7\}$.

Il reste \`{a} v\'{e}rifier la conjecture $(H)$ pour $p\in \{3,4\}$

\begin{enumerate}
\item $p=3.$  D'apr\`{e}s \ref{CS1_ConjH}, la conjecture $(H)$ est vraie dans les cas $n=0$
ou $n\geq E(\frac{1+\sqrt{21}}{2})+1=3$, examinons les deux
autres cas qui restent:
\be
\item $n=1$. Posons $(\bigwedge V,d)=(\bigwedge
(x,y_{1},y_{2},y_{3},y_{4}),d)$. D'apr\`{e}s la minimalit\'{e} de $%
(\bigwedge V,d)$\ on peut supposer que:\newline
$dx=0,dy_{i}=P_{i}(x)(1\leq i\leq 3)$, et $dy_{4}=P_{4}(x)+\alpha $, o\`{u} $\alpha \in 
\mathbb{Q}
\lbrack x]\otimes \bigwedge^{+}\{y_{1},y_{2},y_{3}\}$, et $%
\{P_{1},P_{2},P_{3},P_{4}\}\subset 
\mathbb{Q}
\lbrack x]$\newline
Sans perte de g\'{e}n\'{e}ralit\'{e} on peut poser:\newline
$P_{1}(x)=x^{q},P_{2}(x)=ax^{r}$ et $P_{3}(x)=bx^{s},2\leq q\leq r\leq s.$\textit{\newline
}
Si $a=b=0$, alors l'ensemble $B=\{[y_{2}],[xy_{2}],[y_{3}],[xy_{3}]\}$ est une famille libre dans $H^{impair}(\bigwedge V,d)$, et donc
$$
\dim H^{\star }(\bigwedge V,d)=2\dim H^{impair}(\bigwedge V,d)\geq 8>\dim V=5.
$$

Si $(a,b)\neq (0,0)$, supposons par exemple $a\neq 0$%
. Consid\'{e}rons les \'{e}l\'{e}ments de $H^{impair}(\bigwedge
V,d) $ suivants:
$$\begin{array}{l}
\omega _{1}=[ax^{r-q}y_{1}-y_{2}]\\
\omega_{2}=[bx^{s-q}y_{1}-y_{3}]\\
\omega _{3}=x\omega_{1}=[ax^{r-q+1}y_{1}-xy_{2}]
\end{array}$$

Comme $\{\omega _{1},\omega _{2},\omega _{3}\}$ sont lin%
\'{e}airement ind\'{e}pendants, alors:
$$
\dim H^{\star }(\bigwedge V,d)=2\dim H^{impair}(\bigwedge V,d)\geq 6>\dim V=5.
$$

\item $n=2$. Posons $(\bigwedge V,d)=(\bigwedge
(x_{1},x_{2},y_{1},y_{2},y_{3},y_{4},y_{5}),d)$. Comme ci-dessus on
peut supposer que:
$$\begin{array}{l}
dx_{1}=dx_{2}=0\\
dy_{i}=P_{i}(x_{1},x_{2})\in \mathbb{Q}\lbrack x_{1},x_{2}]\;(i=1,2,3)\\
dy_{j}\in \bigwedge V\;(j=4,5)
\end{array}$$

Posons $B_{2}=\{[x_{1}^{2}],[x_{2}^{2}],[x_{1}x_{2}]\}$ et 
$W_{2}$ le sous espace vectoriel de $H^{pair}(\bigwedge V,d)$%
engendr\'{e} par $B$.

Si $\dim W\geq 1$, alors: 

$$\dim H^{\star }(\bigwedge
V,d)=2\dim H^{pair}(\bigwedge V,d)\geq 8>\dim V=7$$ car $
\mathbb{Q}\{1,[x_{1}],[x_{2}]\}\oplus W_{2}$
est une somme directe dans $%
H^{pair}(\bigwedge V,d)$. 

Supposons maintenant que $W_2=\{0\}$. On peut poser:
$$
\begin{array}{l}
dy_{1} =x_{1}^{2} \\
 dy_{2}=x_{1}x_{2} \\
dy_{3} = x_{2}^{2} \\ 
dy_{4} = P_{4}(x_{1},x_{2})+\alpha \\
 dy_{5} = P_{5}(x_{1},x_{2})+\beta  
\end{array}$$

avec $\alpha ,\beta \in \mathbb{Q}\lbrack x_{1},x_{2}]\otimes \bigwedge^{+}\{y_{1},y_{2},y_{3},y_{4}\}$ qui ne sont pas simultan\'{e}ment nuls (car $X$ n'est pas pur);
supposons par exemple $\alpha \neq 0,\alpha $ s'\'{e}crit:
$$\alpha
=P(x_{1},x_{2})y_{1}y_{2}+Q(x_{1},x_{2})y_{1}y_{3}+R(x_{1},x_{2})y_{2}y_{3}$$
On a alors:
$$
0=d\alpha
=-(x_{2}^{2}Q+x_{1}x_{2}P)y_{1}+(x_{1}^{2}P-x_{2}^{2}R)y_{2}+(x_{1}^{2}Q+x_{1}x_{2}R)y_{3}
$$

ce qui donne le syst\`{e}me:
$$
\begin{array}{lll}
x_{2}^{2}Q+x_{1}x_{2}P & = & 0 \\ 
x_{1}^{2}P-x_{2}^{2}R & = & 0 \\ 
x_{1}^{2}Q+x_{1}x_{2}R & = & 0.%
\end{array}$$

Comme les polyn\^{o}mes $x_{1}^{2}$ et $x_{2}^{2}$ sont premiers entre eux alors il existe $Z\in \mathbb{Q}\lbrack x_{1},x_{2}]$ tel\ que:
$$\begin{array}{l}
P=x_{2}^{2}Z\\Q=-x_{1}x_{2}Z\\
R=x_{1}^{2}Z
\end{array}$$
D'o\`{u}: $\alpha =d(Zy_{1}y_{2}y_{3})$. D'autre part 
$P=d\alpha _{1}$, car $[P]\in W=\{0\}$. On en d\'{e}duit le
cocycle $\gamma =y_{4}-\alpha _{1}-Zy_{1}y_{2}y_{3}$. Consid\'{e}rons les \'{e}l\'{e}ments de $H^{impair}(\bigwedge V,d)$ suivants:
$$\begin{array}{l}
\omega _{1}^{\prime }=[y_{4}-\alpha _{1}-Zy_{1}y_{2}y_{3}]\\
\omega _{2}^{\prime }=x_{1}\omega _{1}^{\prime }\\
\omega_{3}^{\prime }=[x_{2}y_{1}-x_{1}y_{2}]\\
\omega _{4}^{\prime}=[x_{1}y_{3}-x_{2}y_{2}]
\end{array}$$

On v\'{e}rifie facilement que $\omega _{1}^{\prime },\omega
_{2}^{\prime },\omega _{3}^{\prime }$ et $\omega
_{4}^{\prime }$ sont lin\'{e}airement ind\'{e}pendants et donc:
$$\dim H^{\star }(\bigwedge V,d)=2\dim H^{impair}(\bigwedge V,d)\geq 8>\dim V=7.
$$
\ee
\item $p=4.$ On sait que d'après \ref{CS1_ConjH} la conjecture $(H)$ est vraie pour tout $n\leq 2$, ou $n\geq E(\frac{1+\sqrt{33}}{2})+1=4$, il reste donc le cas $n=3$. Posons:\newline
$(\bigwedge V,d)=(\bigwedge
(x_{1},x_{2},x_{3},y_{1},y_{2},y_{3},y_{4},y_{5},y_{6},y_{7}),d)$,
avec $ |x_{i}|\; (1\leq i\leq 3)$ est pair, et $|y_{j}|\; (1\leq j\leq 7)$%
 est impair. Comme ci-dessus $W_{2}$ d\'{e}signe le sous
espace vectoriel de $H^{pair}(\bigwedge V,d)$ engendr\'{e} par%
\newline
$B_{2}=\{e_{ij}=[x_{i}x_{j}]/1\leq i,j\leq 3\}$. 
\\Si $\dim W_{2}\geq
1$, alors:
$$
\dim H^{\star }(\bigwedge V,d)=2\dim H^{pair}(\bigwedge V,d)\geq 10=\dim V.
$$

car $
\mathbb{Q}\{1,[x_{1}],[x_{2}],[x_{3}]\}\oplus W_{2}$ est une somme directe
dans $H^{pair}(\bigwedge V,d)$.
\\Supposons maintenant que
$W_{2}=\{0\}$. On peut \'{e}crire:
$$\begin{array}{ll}
dy_{1} = x_{1}^{2} & dy_{2} = x_{2}^{2} \\
dy_{3} = x_{3}^{2} &dy_{4} =x_{1}x_{2} \\ 
dy_{5} = x_{1}x_{3} & dy_{6} = x_{2}x_{3} \\
dy_{7} = P+\alpha  
\end{array}
$$

o\`{u} $P\in \mathbb{Q}\lbrack x_{1},x_{2},x_{3}]$, et $\alpha \in 
\mathbb{Q}
\lbrack x_{1},x_{2},x_{3}]\otimes
\bigwedge^{+}\{y_{1},y_{2},y_{3},y_{4},y_{5},y_{6}\}.$ Consid\'{e}%
rons les \'{e}l\'{e}ments de $H^{impair}(\bigwedge V,d)$ suivants:
$$
\begin{array}{lllll}
\omega _{1}^{\prime \prime }=[x_{2}y_{1}-x_{1}y_{4}] &,& \omega _{2}^{\prime
\prime }=[x_{3}y_{1}-x_{1}y_{5}]&, & \omega _{3}^{\prime \prime
}=[x_{1}y_{2}-x_{2}y_{4}] \\ 
\omega _{4}^{\prime \prime }=[x_{3}y_{2}-x_{2}y_{6}] &,& \omega _{5}^{\prime
\prime }=[x_{1}y_{3}-x_{3}y_{5}] &,& \omega _{6}^{\prime \prime
}=[x_{2}y_{3}-x_{3}y_{6}].%
\end{array}$$

Il est clair que $\omega _{1}^{\prime \prime }$,  
$\omega
_{2}^{\prime \prime },\omega _{3}^{\prime \prime },\omega _{4}^{\prime \prime },\omega _{5}^{\prime \prime }$ et $\omega _{6}^{\prime \prime }$ sont linèairement
indépendants et donc:

$$\dim H^{\star }(\bigwedge V,d)=2\dim H^{pair}(\bigwedge V,d)\geq 12>\dim V.
$$

\end{enumerate}
\textbf{Troisi\`{e}me cas:} $rk_{0}(X)=-\chi _{\pi }(X)-2=p-2,\;p\geq 2$. Dans ce cas on a:
$$\begin{array}{ll}
\mathit{\dim H}^{\star }\mathit{(}\bigwedge \mathit{V,d)\geq 2}^{p-2}&\text{d'après   théorème }\ref{thm_Hi}\\
\dim H^{\star }(\bigwedge V,d)\geq \dim V,\;\forall n\geq \frac{1+\sqrt{12p-15}}{2}&\text{d'après }(\ref{CS_ConjH})
\end{array}$$
D'o\`{u} :$\dim H^{\star }(\bigwedge V,d)\geq \dim V,\; \forall n\in [0,\frac{2^{p-2}-p}{2}]\cup \lbrack \frac{1+\sqrt{12p-15}}{2}%
,+\infty \lbrack $. Comme ci-dessus cherchons les
valeurs de $p$ pour lesquelles:
$$
\frac{2^{p-2}-p}{2}\geq \frac{1+\sqrt{12p-15}}{2}
$$
c'est \`{a} dire, pour quelles valeurs de $p$ on a: $$A_{p}=2^{p-1}(2^{p-3}-p-1)+p^{2}-10p+16\geq 0$$
On v\'{e}rifie que $A_{p}>0$, pour tout $p\geq 6$. Il
reste \`{a} \'{e}tudier le cas $p\leq 5$. Le cas $p\leq 4$ est d\'{e}j\`{a} trait\'{e} ci-dessus, voyons maintenant le cas $p=5$.

D'après (\ref{CS_ConjH}), la conjecture $(H)$ est v\'{e}rifi\'{e}e pour tout $n\geq E(\frac{1+\sqrt{45}}{2})+1=4$ ou $n\leq E(\frac{2^{3}-5}{2})=1.$

Etudions maintenant les deux cas restants $n=2$ et $n=3$.
\be
\item $n=2$, \'{e}crivons $(\bigwedge
V,d)=(\bigwedge (x_{1},x_{2},y_{1},y_{2},y_{3},y_{4},y_{5},y_{6},y_{7}),d)$ avec $|x_{1}|\leq |x_{2}|$ et $|y_{1}|\leq |y_{2}|\leq |y_{3}|\leq |y_{4}|\leq |y_{5}|\leq |y_{6}|\leq
|y_{7}|.$
Posons pour tout $k\geq 1,\;B_{k}=\{[x_{1}^{i}x_{2}^{j}]/i+j=k\},\;B_{0}=\{1\}$ et $W_{k}$ le sous espace vectoriel de $H^{pair}(\bigwedge V,d)$ engendr\'{e} par $B_{k}$. On a \'{e}videmment $\dim W_{0}=1$ et $\dim W_{1}=2$, trois cas \`{a} distinguer:
\be
\item $\dim (W_{2}+W_{3})\geq 2$,  alors la conjecture $(H)$ est bien v\'{e}rifi\'{e}e, car
$$
\dim H^{\star }(\bigwedge V,d)= 2\dim H^{pair}(\bigwedge V,d)\geq
2\dim (W_{0}\oplus W_{1}\oplus W_{2}+W_{3})\geq 10>\dim V.
$$
\item $\dim (W_{2}+W_{3})=1$, alors $\dim W_{2}=1$ et $W_{k}=\{0\},\;\forall k\geq 3$. Posons:
$$dy_{i}=P_{i}+\alpha _{i}\tq P_{i}\in \mathbb{Q}\lbrack x_{1},x_{2}],\;\alpha _{i}\in
\mathbb{Q}\lbrack x_{1},x_{2}]\otimes \bigwedge{}^{+}(y_{1},...,y_{i-1}).$$
On peut supposer que: $\alpha _{1}=\alpha _{2}=0,\;\alpha
_{3}=Q_{3}y_{1}y_{2}$ et $P_{1},P_{2}$ lin\'{e}airement ind\'{e}pendants dans $\mathbb{Q}\lbrack x_{1}^{2},x_{2}^{2},x_{1}x_{2}]$. Comme $d\alpha _{3}=0$ alors $Q_{3}=0$ et donc $dy_{3}=P_{3}$. Or $W_{3}=\{0\}$, on peut alors prendre $P_{3}\in
\{x_{1}^{3},x_{2}^{3},x_{1}^{2}x_{2},x_{1}x_{2}^{2}\}$. Ecrivons: $\alpha _{4}=Q_{4}y_{1}y_{2}+R_{4}y_{2}y_{3}+T_{4}y_{3}y_{1}$, alors
$0=d^{2}y_{4}=(T_{4}P_{3}-Q_{4}P_{2})y_{1}+(Q_{4}P_{1}-R_{4}P_{3})y_{2}+(R_{4}P_{2}-T_{4}P_{1})y_{3} 
$, d'o\`{u} le syst\`{e}me:
$$
\begin{array}{ccc}
\mathit{T}_{4}\mathit{P}_{3} & \mathit{=} & \mathit{Q}_{4}\mathit{P}_{2} \\ 
\mathit{Q}_{4}\mathit{P}_{1} & \mathit{=} & \mathit{R}_{4}\mathit{P}_{3} \\ 
\mathit{R}_{4}\mathit{P}_{2} & \mathit{=} & \mathit{T}_{4}\mathit{P}_{1}%
\mathit{.}%
\end{array}%
$$

Or $P_{1}$et $P_{2}$ sont premiers entre eux,
donc il existe $Q\in \mathbb{Q}
\lbrack x_{1},x_{2}]$ tel que: $R_{4}=QP_{1}$ et $T_{4}=QP_{2}$, ainsi $Q_{4}=QP_{3}$ , d'o\`{u} $dy_{4}=P_{4}+d(Qy_{1}y_{2}y_{3})$. D'autre part on a: $[P_{4}]=0$, ce qui implique que $P_{4}=d(P_{1}^{\prime }y_{1}+P_{2}^{\prime
}y_{2}+P_{3}^{\prime }y_{3})$. Par cons\'{e}quent:
$$
d(y_{4}-P_{1}^{\prime }y_{1}+P_{2}^{\prime }y_{2}+P_{3}^{\prime
}y_{3}-Qy_{1}y_{2}y_{3})=0.
$$

Consid\'{e}rons $y=y_{4}-P_{1}^{\prime }y_{1}+P_{2}^{\prime
}y_{2}+P_{3}^{\prime }y_{3}-Qy_{1}y_{2}y_{3},\;\omega =[P]$ un g\'{e}n\'{e}rateur de $W_{2}$, et $\mu $ la classe
fondamentale de $H^{\star }(\bigwedge V,d)$ (ie: un g\'{e}n\'{e}%
rateur de $H^{fd(X)}(\bigwedge V,d)$. Alors les \'{e}l\'{e}ments
suivants: $[y],[x_{1}y],[x_{2}y],[Py]$ et $\mu $ sont lin\'{e}airement ind\'{e}pendants dans $H^{impair}(\bigwedge V,d)$, d'o\`{u} la conjecture $(H)$ est vérifiée.
\item $W_{3}=\{0\}$, on peut alors \'{e}crire:
$$
\begin{array}{l}
\mathit{dy}_{1} =\mathit{x}_{1}^{2} \\
\mathit{dy}_{2} = \mathit{x}_{1}\mathit{x}_{2} \\ 
\mathit{dy}_{3} = \mathit{x}_{2}^{2} \\
\mathit{dy}_{4} = \mathit{P}_{4}\mathit{+d(Qy}_{1}\mathit{y}_{2}\mathit{y}_{3}).
\end{array}
$$
On conclut comme ci-dessus.
\ee
\item $n=3$, \'{e}crivons $(\bigwedge V,d)=(\bigwedge (x_{1},x_{2},x_{3},y_{1},y_{2},y_{3},y_{4},y_{5},y_{6},y_{7},y_{8}),d)$
, avec
$|x_{1}|\leq |x_{2}|\leq |x_{3}|$, et $|y_{1}|\leq |y_{2}|\leq
|y_{3}|\leq |y_{4}|\leq |y_{5}|\leq |y_{6}|\leq |y_{7}|\leq |y_{8}|.$

Posons $B_{k}^{\prime }=\{[x_{1}^{i}x_{2}^{j}x_{3}^{l}]/i+j+l=k\},B_{0}^{\prime }=\{1\}$ et $W_{k}^{\prime }$ le sous
espace vectoriel de $H^{pair}(\bigwedge V,d)$ engendr\'{e} par $B_{k}$. 
\\On a $\dim (W_{0}^{\prime }\oplus W_{1}^{\prime })=4$.
\be
\item Si $\dim (W_{2}^{\prime }+W_{3}^{\prime })\geq 2$, alors: 
$$\dim H^{pair}(\bigwedge V,d)\geq 2[\dim (W_{0}^{\prime }\oplus W_{1}^{\prime
})+\dim (W_{2}^{\prime }+W_{3}^{\prime })]\geq 12>\dim V$$
\item Si $\dim (W_{2}^{\prime }+W_{3}^{\prime })=1$, alors $\dim W_{2}^{\prime }=1$ et $W_{k}^{\prime }=\{0\},\forall
k\geq 3$. Posons:
$$dy_{i}=P_{i}+\alpha _{i}\tq P_{i}\in 
\mathbb{Q},\;\lbrack x_{1},x_{2},x_{3}],\alpha _{i}\in 
\mathbb{Q}
\lbrack x_{1},x_{2},x_{3}]\otimes \bigwedge\;{}^{+}(y_{1},...,y_{i-1}).$$

On peut alors supposer d'apr\`{e}s les hypoth\`{e}ses que:
\\$dy_{i}=P_{i},\;(1\leq i\leq 5)$ sont lin\'{e}airement ind\'{e}pendants dans $
\mathbb{Q}
\{x_{1}^{2},x_{2}^{2},x_{3}^{2},x_{1}x_{2},x_{1}x_{3},x_{2}x_{3}\}$ et $dy_{6}=P_{6}\in \{x_{1}^{i}x_{2}^{j}x_{3}^{l}/i+j+l=3\}.$ Soit $\alpha =[P]$ un g\'{e}n\'{e}rateur de $W_{2}^{\prime }$, on peut prendre $P=x_{i}x_{j}$, o\`{u} $1\leq i\leq j\leq 3$. On peut se restreindre aux deux cas: $(i,j)=(1,1)$ ou $(i,j)=(1,2)$.
\be
\item $P=x_{1}^{2}$, on peut alors \'{e}crire:
$$
\begin{array}{llllll}
\mathit{P}_{1} & \mathit{=} & \mathit{a}_{1}\mathit{x}_{1}^{2}\mathit{+x}%
_{2}^{2} & \mathit{P}_{2} & \mathit{=} & \mathit{a}_{2}\mathit{x}_{1}^{2}%
\mathit{+x}_{3}^{2} \\ 
\mathit{P}_{3} & \mathit{=} & \mathit{a}_{3}\mathit{x}_{1}^{2}\mathit{+x}_{1}%
\mathit{x}_{2} & \mathit{P}_{4} & \mathit{=} & \mathit{a}_{4}\mathit{x}%
_{1}^{2}\mathit{+x}_{1}\mathit{x}_{3} \\ 
\mathit{P}_{5} & \mathit{=} & \mathit{a}_{5}\mathit{x}_{1}^{2}\mathit{+x}_{2}%
\mathit{x}_{3} & \mathit{P}_{6} & \mathit{=} & \mathit{x}_{1}^{3}\mathit{.}%
\end{array}%
$$

Par un simple calcul on v\'{e}rifie que les \'{e}l\'{e}ments
suivants sont des cocycles
$$
\begin{array}{lll}
\alpha _{1} & = & 
x_{2}y_{2}+a_{2}x_{1}y_{3}+a_{5}x_{1}y_{4}-x_{3}y_{5}+(a_{2}a_{3}-a_{4}a_{5})y_{6}
\\ \alpha _{2} & = & 
x_{3}y_{3}-a_{3}x_{1}y_{4}-x_{1}y_{5}-(a_{3}a_{4}+a_{5})y_{6} \\ 
\alpha _{3} & = & 
x_{3}y_{1}+a_{5}x_{1}y_{3}-a_{1}x_{1}y_{4}-x_{2}y_{5}+(a_{1}a_{4}-a_{3}a_{5})y_{6}
\\ \alpha _{4} & = & x_{2}y_{1}+(a_{3}x_{1}-x_{2})y_{3}-(a_{1}+a_{3}^{2})y_{6}
\\ 
\alpha _{5} & = & 
a_{4}x_{1}y_{3}-x_{2}y_{4}+x_{1}y_{5}-(a_{3}a_{4}+a_{5})y_{6} \\
 \alpha _{6}
& = & x_{1}y_{2}+(a_{4}x_{1}-x_{3})y_{4}-(a_{2}+a_{4}^{2})y_{6}.%
\end{array}%
$$

Consid\'{e}rons la matrice:
$$M=\left( 
\begin{array}{cccccc}
0 & 0 & 0 & 0 & 0 & 0 \\ 
0 & 0 & 0 & 1 & 0 & 0 \\ 
0 & 0 & 1 & 0 & 0 & 0 \\ 
0 & 0 & 0 & 0 & 0 & 1 \\ 
1 & 0 & 0 & 0 & 0 & 0 \\ 
0 & 0 & 0 & 0 & 0 & 0 \\ 
a_{2} & 0 & a_{5} & a_{3} & a_{4} & 0 \\ 
0 & 0 & 0 & -1 & 0 & 0 \\ 
0 & 1 & 0 & 0 & 0 & 0 \\ 
a_{5} & -a_{3} & -a_{1} & 0 & 0 & a_{4} \\ 
0 & 0 & 0 & 0 & -1 & 0 \\ 
0 & 0 & 0 & 0 & 0 & -1 \\ 
0 & -1 & 0 & 0 & 1 & 0 \\ 
0 & 0 & -1 & 0 & 0 & 0 \\ 
-1 & 0 & 0 & 0 & 0 & 0 \\ 
(a_{2}a_{3}-a_{4}a_{5}) & -(a_{3}a_{4}+a_{5}) & (a_{1}a_{4}-a_{3}a_{5}) & 
-(a_{1}+a_{3}^{2}) & -(a_{3}a_{4}+a_{5}) & -(a_{2}+a_{4}^{2})%
\end{array}%
\right) $$

Soient $E$ le $\mathbb{Q}$-espace vectoriel engendr\'{e} par l'ensemble \\$B=\{x_{1}y_{1},x_{2}y_{1},x_{3}y_{1},x_{1}y_{2},x_{2}y_{2},x_{3}y_{2},x_{1}y_{3},x_{2}y_{3},x_{3}y_{3},\\x_{1}y_{4},x_{2}y_{4},x_{3}y_{4},x_{1}y_{5},x_{2}y_{5},x_{3}y_{5},y_{6}\} 
$ et $\{e_{1},e_{2},e_{3},e_{4},e_{5},e_{6}\}$ la base
canonique de $\R^{6}$. Alors on a:
$\alpha _{i}=Me_{i},\;\forall i,\;\;1\leq i\leq 6$ (les
composantes sont dans la base $B$)

La matrice $M$ \'{e}tant de rang 6, donc $[\alpha _{1}],[\alpha _{2}],[\alpha _{3}],[\alpha
_{4}] ,[\alpha _{5}]$, et $[\alpha _{6}]$ sont
lin\'{e}airement ind\'{e}pendants dans $H^{impair}(\bigwedge V,d)$, d'o\`{u} la conjecture $(H)$ est réalisée.
\item $P=x_{1}x_{2}$, posons alors:
$$\begin{array}{llllll}
P_{1} & = & a_{1}x_{1}x_{2}+x_{1}^{2} & P_{2} & = & a_{2}x_{1}x_{2}+x_{2}^{2}
\\ 
P_{3} & = & a_{3}x_{1}x_{2}+x_{3}^{2} & P_{4} & = & 
a_{4}x_{1}x_{2}+x_{1}x_{3} \\ 
P_{5} & = & a_{5}x_{1}x_{2}+x_{2}x_{3} & P_{6} & = & x_{1}^{2}x_{2}.%
\end{array}$$

On montre de la m\^{e}me mani\`{e}re, que les \'{e}l\'{e}ments
suivants sont lin\'{e}airement ind\'{e}pendants dans $H^{impair}(\bigwedge
V,d)$, 
$$\begin{array}{l}
\left[x_{2}y_{1}-a_{1}x_{1}y_{2}+(a_{1}a_{2}-1)y_{6}\right]\\
\left[x_{3}y_{1}-x_{1}y_{4}-a_{1}x_{1}y_{5}+(a_{1}a_{5}+a_{4})y_{6}\right]\\
\left[x_{2}y_{4}-a_{4}x_{1}y_{2}-x_{1}y_{5}+(a_{2}a_{4}+a_{5})y_{6}\right]\\
\left[x_{3}y_{4}-x_{1}y_{3}-a_{4}x_{1}y_{5}+(a_{4}a_{5}+a_{3})y_{6}\right]\\
\left[(a_{5}x_{1}+x_{3})y_{2}-(a_{2}x_{1}+x_{2})y_{5}-2a_{2}a_{5}y_{6}\right]\\
\left[a_{3}x_{1}y_{2}-x_{2}y_{3}+(x_{3}-a_{5}x_{1})y_{5}-(a_{2}a_{3}-a_{5}^{2})y_{6}\right]
\end{array}$$

D'o\`{u}  la conjecture $(H)$
\ee
\item $W_{2}^{\prime }=\{0\}$, dans ce cas on peut poser:
$$\begin{array}{lllllllll}
dy_{1} & = & x_{1}^{2} & dy_{2} & = & x_{2}^{2} & dy_{3} & = & x_{3}^{2} \\ 
dy_{4} & = & x_{1}x_{2} & dy_{5} & = & x_{1}x_{3} & dy_{6} & = & x_{2}x_{3}.%
\end{array}$$

Les \'{e}l\'{e}ments suivants sont lin\'{e}airement ind\'{e}pendants
dans $H^{impair}(\bigwedge V,d)$:\newline
$[x_{2}y_{1}-x_{1}y_{4}]$\textit{, }$[x_{1}y_{2}-x_{2}y_{4}]$\textit{, }$%
[x_{3}y_{1}-x_{1}y_{4}]$\textit{, }$[x_{1}y_{3}-x_{3}y_{5}]$\textit{, }$%
[x_{3}y_{2}-x_{2}y_{6}]$ et $[x_{2}y_{3}-x_{3}y_{6}]$, ce
qui ach\`{e}ve la d\'{e}monstration. \hfill$\Box $
\ee
\ee

\vskip 4mm
\textbf{Preuve Théorème E}.
Posons $fd(X)=N$, on a d'après théorème \ref{thm_FH}  que $\ds\sum_{i=1}^n|x_i|\ioe N$ et $\ds\sum_{i=1}^{n+p}|y_i|\ioe 2N-1$, or $|x_i|\soe 2$ et $|y_i|\soe 3$, d'où $n \ioe\dfrac N2,n+p\ioe \dfrac {2N-1}3$ et $\dim V=2n+p\ioe \dfrac {7N-2}6$, d'autre part on sait théorème \ref{thm_Hi} que  $\dim H^*(X,\Q)\soe 2^{rk_0(X)}$ or $rk_0(X)=N-{\rm codim}(X)\soe N-6$, donc $\dim H^*(X,\Q)\soe 2^{N-6}$, on cherche alors les entiers vérifiant $2^{N-6}\soe \dfrac {7N-2}6$ ou bien $ 3. 2^N-224N+64\soe 0$

L'étude de la fonction $ f(N)=3. 2^N-224N+64$ montre qu'elle est croissante sur l'intervalle $\left[\dfrac{\frac{224}{\ln 8}}{\ln 2},+\infty\right[$ or $\dfrac{\frac{224}{\ln 8}}{\ln 2}\ioe 10$, donc $\qlq N\soe 10,\;f(N)\soe f(10)=896$, alors que pour $N\ioe 10$, la conjecture (H) est vrai d'après \ref{CS3_ConjH}  \hfill $\Box$

\subsection{Cas symplectique} (avec condition)

\vskip 3mm
\textbf{Preuve du Théorème F}.
On rappelle que $X$ les espaces symplectiques vérifient les propriétés suivantes \cite{AP}:
$$\begin{array}{l}
fd(X)=2m\\
\exists w\in H^2(X,\Q)\tq H^{2m}(X,\Q)\cong \Q w^m\\
\text{le cup-produit }w^k:H^{m-k}(X,\Q)\vers H^{m+k}(X,\Q)\\
\text{est un isomorphisme}\\
\end{array}$$

Ainsi les nombres de Betti $b_{2i}=\dim H^{2i}(X,\Q), 0\ioe i\ioe m$ sont tous non nuls, donc $\dim H^{\text{pair}}(X,\Q)\soe m+1$, or $\chi_c=0$ car le cas contraire ($\chi_\pi=0$) est déjà résolu dans le cas pur, donc 
$$\begin{array}[t]{lll}
\dim H^*(X,\Q)&=2\dim H^{\text{pair}}(X,\Q)&\\
&\soe 2m+2&\\
&>2m&\\
&=fd(X)\soe \dim V&\text{d'après théorème \ref{thm_FH}}
\end{array}$$

\hfill $\Box$

\subsection{Cas cosymplectique} (sans condition)

\vskip 3mm

\textbf{Preuve du Théorème G}.
Si $X=M^{2n+1}$ est une variété co-symplectique, d'aprés \cite{BG}  les nombres de Betti sont tous non nuls,  plus encore, dans\cite{CLM} on a: 
$b_0=b_{2n+1}\ioe b_1=b_{2n}\ioe\cdots\ioe b_n=b_{n+1}$ donc $\dim H^*(\wedge V,\Q)=\ds\sum_{i=0}^{2n+1}b_i\soe 2n+2>2n+1=fd(X)\soe \dim V$. \hfill $\Box$

\vskip 4mm
\begin{remark}
D'après la définition d'une variété cosymplectique  il existe une I-forme
fermée $\eta $ et donc une classe dans le $H^1(X)$ qui n'est pas nulle.

Tous les exemples connus de variétés
cosymplectiques sont NON 1-connexes. Le $\pi_1$ est celui d'un
tore de dimension
$\soe 1$ ou un groupe nilpotent $(X = G/\Gamma )$ est  une nilvariété.

Cependant  il s'agit d'espaces simples ($\pi_1$ abélien opérant trivialement
sur les $\pi_n$) ou bien d'espaces nilpotents ($\pi_1$ nilpotent opérant de
manière nilpotente sur les $\pi_n$).  

Ces espaces admettent des modèles
minimaux de Sullivan (non 1 connexe) qui vérifie la propriété 
 $fd (X ) \soe dim V$
\end{remark}

\vskip 10mm
\end{document}